# NECESSARY AND SUFFICIENT CONDITION FOR ASYMPTOTIC NORMALITY OF STANDARDIZED SAMPLE MEANS

BY RAJESHWARI MAJUMDAR* AND SUMAN MAJUMDAR

*University of Connecticut and University of Connecticut*

The double sequence of standardized sample means constructed from an infinite sequence of square integrable independent random vectors in the plane with identically distributed coordinates is jointly asymptotically Normal if and only if the Cesaro means of the sequence of cross sample correlation coefficients converges to 0.

**1. Introduction and Results.** We investigate the joint asymptotic Normality of the standardized sample means, based on two random samples from two distributions, in this paper. Hereinafter, iid will abbreviate independent and identically distributed. Let

$$\{X_{1,j} : j \geq 1\} - \text{iid sequence of random variables, mean } \mu_1 \text{ and variance } \sigma_1^2 \quad (1)$$

$$\{X_{2,j} : j \geq 1\} - \text{iid sequence of random variables, mean } \mu_2 \text{ and variance } \sigma_2^2. \quad (2)$$

In what follows, we always assume that the index $i$ runs from to 2. Let

$$Y_{n_i}^{(i)} = \frac{\sqrt{n_i}(\bar{X}_i - \mu_i)}{\sigma_i} = \frac{1}{\sqrt{n_i}} \sum_{j=1}^{n_i} \left( \frac{X_{ij} - \mu_i}{\sigma_i} \right) \quad (3)$$

denote the standardized sample mean from the $i^{\text{th}}$ sample. Let $\Phi$ denote the standard Normal distribution on the line and $\Phi \times \Phi$ the product measure on the plane, the bivariate standard Normal distribution. Our objective is to obtain necessary and sufficient conditions under which, as $n_1, n_2 \to \infty$, the double sequence of random vectors

$$Y_{n_1,n_2} = \left( Y_{n_1}^{(1)}, Y_{n_2}^{(2)} \right) \text{ converges in distribution to } \Phi \times \Phi. \quad (4)$$

Let $\psi_{(i)}$ denote the characteristic function (CF, hereinafter) of the standardized $X_{i,j}$, which, by (1) and (2), does not depend on $j$, and $f_{n_i}^{(i)}$ that of $Y_{n_i}^{(i)}$; note that

$$f_{n_i}^{(i)}(u) = \left[ \psi_{(i)} \left( \frac{u}{\sqrt{n_i}} \right) \right]^{n_i}. \quad (5)$$

Since, for every $u \in \Re$, $x \mapsto \exp(iux)$ is bounded and continuous, where $i^2 = -1$, by the Central Limit Theorem (CLT, hereinafter) [Dudley (1989, Theorem 9.5.6)] and the definition of convergence in distribution, with $\zeta$ denoting the CF of $\Phi$ so that $\zeta(u) = \exp(-u^2/2)$,

$$\lim_{n_i \to \infty} f_{n_i}^{(i)}(u) = \zeta(u). \tag{6}$$

Let $\psi_{n_1, n_2}$ denote the CF of $Y_{n_1, n_2}$ and $\psi$ that of $\Phi \times \Phi$, so that

$$\psi((s, t)) = \exp\left(-\left(s^2 + t^2\right)/2\right) = \zeta(s)\zeta(t).$$

Assume that

the sequence $\{X_{1,j} : j \geq 1\}$ is independent of the sequence $\{X_{2,j} : j \geq 1\}$. \tag{7}

By (7), $\psi_{n_1, n_2}((s, t)) = f_{n_1}^{(1)}(s) f_{n_2}^{(2)}(t)$. Since $\psi((s, t)) = \zeta(s)\zeta(t)$,

$$\left|\psi_{n_1, n_2}((s, t)) - \psi((s, t))\right|$$
$$\leq \left|f_{n_1}^{(1)}(s) f_{n_2}^{(2)}(t) - f_{n_1}^{(1)}(s)\zeta(t)\right| + \left|f_{n_1}^{(1)}(s)\zeta(t) - \zeta(s)\zeta(t)\right|;$$

since every CF is bounded, by (6), for every $(s, t) \in \Re^2$,

$$\lim_{n_1, n_2 \to \infty} \psi_{n_1, n_2}((s, t)) = \psi((s, t)). \tag{8}$$

If $n_1 = n_2$, then (4) follows from (8) by the Levy Continuity Theorem (LCT, hereinafter) [Dudley (1989, Theorem 9.8.2)].

Does (8) imply (4) even without the restriction $n_1 = n_2$? The answer is yes, but to substantiate that assertion we have to come up with an appropriate LCT. We do that by: restating (4) in terms of weak convergence of induced measures, recognizing that a double sequence is a net, and formulating the LCT for nets of random vectors. The edifice thus constructed is used in establishing the main result, Theorem 1, which weakens the collection of assumptions (1), (2), and (7) that is sufficient for (4) to a collection of assumptions that includes a constituent which, in the presence of the remaining constituents, is necessary and sufficient for (4).

For a separable metric space $\mathcal{S}$, let $\mathcal{B}(\mathcal{S})$ denote the Borel $\sigma$-algebra of $\mathcal{S}$ and $\mathcal{M}(\mathcal{S})$ the set of probability measures on $\mathcal{B}(\mathcal{S})$. Endowed with the topology of weak convergence, $\mathcal{M}(\mathcal{S})$ is metrizable as a separable metric space [Parthasarathy (1967, Theorem II.6.2)]. Since all the random elements under consideration are Borel measurable, convergence in distribution is equivalent to weak convergence of the induced probability measures [van der Vaart and Wellner (1996, page 18)].

Let $\mathfrak{N}$ denote the set of natural numbers; then $\mathfrak{N} \times \mathfrak{N}$ is a directed set under the partial ordering $\succeq$ defined by

$$(n_1, n_2) \succeq (m_1, m_2) \Leftrightarrow (n_1 \geq m_1 \text{ and } n_2 \geq m_2).$$

A double sequence $\{x_{n_1, n_2} : n_1, n_2 \geq 1\}$ taking values in $\mathcal{S}$ converges to $x \in S$ as $n_1, n_2 \to \infty$ if and only if the corresponding net $\{x_\alpha : \alpha \in \mathfrak{N} \times \mathfrak{N}\}$ converges to $x$.

Remark 1 formulates the LCT for nets of random vectors.

**Remark 1** Example 1.3.5 of van der Vaart and Wellner (1996) assert, without a proof, the LCT for any $\Re^k$-valued net of random vectors. They remark that a proof can be based on Prohorov's theorem [van der Vaart and Wellner (1996, Theorem 1.3.9)] and the uniqueness of CFs [Dudley (1989, Theorem 9.5.1)]. The proof of the LCT in the case of sequences first establishes that pointwise convergence of CFs implies uniform tightness of the underlying measures (using the Dominated Convergence Theorem) and then uses Prohorov's theorem and the uniqueness of CFs. However, it is not immediately clear how the Dominated Convergence Theorem and consequently, the proof of pointwise convergence of CFs implying uniform tightness in the case of sequences, generalizes to the case of netss. As such, we will use the LCT for nets only after verifying uniform tightness of the underlying measures.                                          //

Remark 2 completes the process of establishing that the collection of assumptions (1), (2), and (7) is sufficient for (4), by verifying the appropriate uniform tightness that would allow us to conclude (4) from (8). For $k \in \mathfrak{N}$, let $\overline{X}_{i,k}$ denote $\sum_{j=1}^{k} X_{i,j}/k$. Let $n_i : \mathfrak{N} \times \mathfrak{N} \mapsto \mathfrak{N}$ denote the order preserving and cofinal map that maps $\alpha$ to its $i^{\text{th}}$ coordinate $n_i(\alpha)$.

**Remark 2** For $\alpha \in \mathfrak{N} \times \mathfrak{N}$, let $Y_\alpha^{(i)}$ denote $Y_{n_i(\alpha)}^{(i)}$ defined in (3) and $H_\alpha \in \mathcal{M}(\Re^2)$ the measure induced by $Y_\alpha = \left( Y_\alpha^{(1)}, Y_\alpha^{(2)} \right)$. Then (4) is equivalent to (in $\mathcal{M}(\Re^2)$)

$$\lim_\alpha H_\alpha = \Phi \times \Phi. \tag{9}$$

Recall that $\psi_\alpha = \psi_{n_1(\alpha), n_2(\alpha)}$ is the CF of $Y_\alpha$, so that the assertion of (8) can be restated as

$$\lim_\alpha \psi_\alpha((s,t)) = \psi((s,t)) \tag{10}$$

for every $(s,t) \in \Re^2$. Once we show

$$\{H_\alpha : \alpha \in \mathfrak{N} \times \mathfrak{N}\} \text{ is uniformly tight,} \tag{11}$$

by Remark 1, the LCT for nets can be used to conclude that (10) implies (9). Note that, since the kernel of the CF is bounded and continuous, (9) implies (10), so that post successful verification of (11), (9) and (10) are equivalent.

For $k \in \mathfrak{N}$ and $i = 1, 2$, let $P_k^{(i)} \in \mathcal{M}(\Re)$ denote the measure induced by $\sqrt{k}(\overline{X}_{i,k} - \mu_i)/\sigma_i$, which, under (1) or (2), converges to $\Phi$ in $\mathcal{M}(\Re)$ by the CLT. By Proposition 9.3.4 of Dudley (1989), $\left\{ P_k^{(i)} : k \in \mathfrak{N} \right\}$ is uniformly tight. Consequently, $\left\{ P_{n_i(\alpha)}^{(i)} : \alpha \in \mathfrak{N} \times \mathfrak{N} \right\}$, being contained in $\left\{ P_k^{(i)} : k \in \mathfrak{N} \right\}$, is uniformly tight as well. By Tychonoff's theorem and Bonneferroni's inequality, (11) follows.

Note that (11) is obtained without (7), that is, by using only (1) and (2).                                          //

We now start the process of stating Theorem 1. The triplet of assumptions (1), (2), and (7) is equivalent to the pair of assumptions

$$\{(X_{1,j}, X_{2,j}) : j \geq 1\} \ - \ \text{iid sequence of random vectors} \tag{12}$$

$$\text{for every } j \geq 1, \, X_{1,j} \text{ and } X_{2,j} \text{ are independent.} \tag{13}$$

We first weaken the assumption in (12) to

$$\{(X_{1,j}, X_{2,j}) : j \geq 1\} \text{ is an independent sequence of random vectors.} \tag{14}$$

Let us introduce some notations here. Let $\xi : \mathfrak{N} \mapsto \mathfrak{N} \times \mathfrak{N}$ denote the order preserving and cofinal map given by $\xi(k) = (k, k)$. For $\alpha \in \mathfrak{N} \times \mathfrak{N}$, let

$$n_\alpha = \min(n_1(\alpha), n_2(\alpha)) \quad \text{and} \quad m_\alpha = \sqrt{n_1(\alpha) n_2(\alpha)}.$$

Let

$$\rho_{jk} \text{ denote the correlation coefficient between } X_{1,j} \text{ and } X_{2,k}, \tag{15}$$

and for $\alpha \in \mathfrak{N} \times \mathfrak{N}$, let

$$\overline{\rho}_\alpha = \frac{1}{m_\alpha} \sum_{j=1}^{n_\alpha} \rho_{jj}. \tag{16}$$

We now weaken the assumption of independence in (13) to either of

$$\lim_\alpha \overline{\rho}_\alpha = 0 \tag{17}$$

$$\lim_k \overline{\rho}_{\xi(k)} = 0. \tag{18}$$

Next, we introduce the assumption

$$\lim_\alpha \widehat{F}_\alpha = \Phi, \tag{19}$$

where $\widehat{F}_\alpha \in \mathcal{M}(\Re)$ is the measure induced on $\mathcal{B}(\Re)$ by

$$\widehat{W}_\alpha = \frac{\left(\overline{X}_{1,n_1(\alpha)} - \mu_1\right) - \left(\overline{X}_{2,n_2(\alpha)} - \mu_2\right)}{\sqrt{\frac{\sigma_1^2}{n_1(\alpha)} + \frac{\sigma_2^2}{n_2(\alpha)}}} = \langle Y_\alpha, v_\alpha \rangle, \tag{20}$$

$v_\alpha = \left(v_\alpha^{(1)}, \, -v_\alpha^{(2)}\right) \in \Re^2$, and

$$v_\alpha^{(i)} = \frac{\sigma_i}{\sqrt{e_\alpha^{1-i} \sigma_1^2 + e_\alpha^{2-i} \sigma_2^2}} \in [0, 1]. \tag{21}$$

**Theorem 1** If (1), (2), and (14) hold, then (17), (18), (9), and (19) are equivalent.

**Remark 3** The collection of assumptions (1), (2), (14), and (18) is weaker than the pair of assumptions (12) and (13). To see that, consider a pair of dependent but uncorrelated random variables and a sequence of iid copies of the resulting random vector.                    //

We present the proof of Theorem 1, an important corollary to it [Proposition 1], and the auxiliary results used in the proof of Theorem 1 in Section 2. The technical results used in all the proofs of this section are assembled in two Appendices, A.1 and A.2.

**2. Proofs.** We prove Theorem 1 by showing

$$(17) \Rightarrow (18) \Rightarrow (9) \Rightarrow (19) \Rightarrow (17). \tag{22}$$

By Lemma A.1, $(17) \Rightarrow (18)$.

The proofs of the remaining assertions in (22) make critical use of the compactness of $[0, \infty]$, the one-point compactification of $[0, \infty)$ [Dudley (1989, Theorem 2.8.1)]. Let

$$e_\alpha = n_1(\alpha)/n_2(\alpha).$$

Since every net taking values in a compact set has a convergent subnet [Lemma A.3], every subnet $\{e_{\phi(\beta)} : \beta \in \mathcal{F}\}$ of $\{e_\alpha : \alpha \in \mathfrak{N} \times \mathfrak{N}\}$ has a further subnet $\{e_{\phi(\varphi(\delta))} : \delta \in \mathfrak{D}\}$ such that

$$\lim_\delta e_{\phi(\varphi(\delta))} = \kappa \in [0, \infty]. \tag{23}$$

For subsequent use, let $(\Phi \times \Phi)_\rho \in \mathcal{M}(\Re^2)$ denote the Normal distribution with means 0, variances 1, and correlation coefficient $\rho \in [-1, 1]$, so that $(\Phi \times \Phi)_0 = \Phi \times \Phi$, and $\mathcal{N}_\theta \in \mathcal{M}(\Re)$ the Normal distribution with mean 0 and variance $\theta \geq 0$, so that $\mathcal{N}_0$ is the point mass at 0 and $\mathcal{N}_1 = \Phi$.

<u>Proof of</u> **(18) $\Rightarrow$ (9)**: Recall from Remark 2 that, under (1) and (2), (9) and (10) are equivalent. Lemma 1 is a key step in proving (18) implies (10).

**Lemma 1** Let $\overline{\rho}_\alpha$ be as in (16). If (14) holds, then, for every subsequence $\{H_{\xi(k(r))} : r \in \mathfrak{N}\}$ of $\{H_{\xi(k)} : k \in \mathfrak{N}\}$, there exists a further subsequence $\{H_{\xi(k(r(m)))} : m \in \mathfrak{N}\}$ such that

$$\lim_m \overline{\rho}_{\xi(k(r(m)))} = \rho \in [-1, 1] \tag{24}$$

and

$$\lim_m H_{\xi(k(r(m)))} = (\Phi \times \Phi)_\rho. \tag{25}$$

The proof of Lemma 1 is presented at the end of the section.

By Lemma 1, (18) implies

$$\lim_k H_{\xi(k)} = \Phi \times \Phi, \tag{26}$$

reducing the proof to showing (26) $\Rightarrow$ (10). Clearly, (26) implies, for every $(u, v) \in \Re^2$,

$$\lim_k \psi_{\xi(k)}((u, v)) = \psi((u, v)). \tag{27}$$

Fix $(s, t) \in \Re^2$ arbitrarily. By Lemma A.2, (10) follows if given any subnet $\{\psi_{\phi(\beta)}((s, t)) : \beta \in \mathcal{F}\}$ of $\{\psi_\alpha((s, t)) : \alpha \in \mathfrak{N} \times \mathfrak{N}\}$, we can find a further subnet $\{\psi_{\phi(\varphi(\delta))}((s, t)) : \delta \in \mathfrak{D}\}$ such that $\psi_{\phi(\varphi(\delta))}((s, t))$ converges to $\psi((s, t))$.

For $\alpha \in \mathfrak{N} \times \mathfrak{N}$, define

$$J_\alpha^{1,2} = \begin{cases} 1 & \text{if } n_1(\alpha) > n_2(\alpha) \\ 0 & \text{otherwise} \end{cases} \text{ and } J_\alpha^{2,1} = \begin{cases} 1 & \text{if } n_1(\alpha) < n_2(\alpha) \\ 0 & \text{otherwise,} \end{cases}$$

so that

$$1 - \left( J_\alpha^{1,2} + J_\alpha^{2,1} \right) = \begin{cases} 1 & \text{if } n_1(\alpha) = n_2(\alpha) \\ 0 & \text{otherwise;} \end{cases}$$

note that $\{J_\alpha^{1,2} : \alpha \in \mathfrak{N} \times \mathfrak{N}\}$ and $\{J_\alpha^{2,1} : \alpha \in \mathfrak{N} \times \mathfrak{N}\}$ are two nets in $\{0, 1\}$. A straightforward algebraic calculation, using (14) and (5), leads to the decomposition

$$\psi_\alpha((s, t)) = \psi_\alpha^{(1)}((s, t)) + \psi_\alpha^{(2)}((s, t)) + \psi_\alpha^{(3)}((s, t)), \tag{28}$$

where

$$\begin{aligned} \psi_\alpha^{(1)}((s, t)) &= \psi_{\xi(n_1(\alpha))}((s, t)) \left( 1 - \left( J_\alpha^{1,2} + J_\alpha^{2,1} \right) \right) \\ \psi_\alpha^{(2)}((s, t)) &= \psi_{\xi(n_2(\alpha))} \left( \left( e_\alpha^{-\frac{1}{2}} s, t \right) \right) \left( \psi_{(1)} \left( \frac{s}{\sqrt{n_1(\alpha)}} \right) \right)^{n_1(\alpha) - n_\alpha} \left( J_\alpha^{1,2} \right) \\ \psi_\alpha^{(3)}((s, t)) &= \psi_{\xi(n_1(\alpha))} \left( \left( s, e_\alpha^{\frac{1}{2}} t \right) \right) \left( \psi_{(2)} \left( \frac{t}{\sqrt{n_2(\alpha)}} \right) \right)^{n_2(\alpha) - n_\alpha} \left( J_\alpha^{2,1} \right). \end{aligned} \tag{29}$$

By Lemma A.3, given an arbitrary subnet $\{\psi_{\phi(\beta)} : \beta \in \mathcal{F}\}$ of $\{\psi_\alpha : \alpha \in \mathfrak{N} \times \mathfrak{N}\}$, we can find a further subnet $\{e_{\phi(\varphi(\delta))} : \delta \in \mathfrak{D}\}$ of $\{e_{\phi(\beta)} : \beta \in \mathcal{F}\}$, $\{J_{\phi(\varphi(\delta))}^{1,2} : \delta \in \mathfrak{D}\}$ of $\{J_{\phi(\beta)}^{1,2} : \beta \in \mathcal{F}\}$, and $\{J_{\phi(\varphi(\delta))}^{2,1} : \delta \in \mathfrak{D}\}$ of $\{J_{\phi(\beta)}^{2,1} : \beta \in \mathcal{F}\}$ such that

$$\lim_\delta J_{\phi(\varphi(\delta))}^{1,2} = J^{1,2} \in \{0, 1\} \text{ and } \lim_\delta J_{\phi(\varphi(\delta))}^{2,1} = J^{2,1} \in \{0, 1\} \text{ exist,} \tag{30}$$

and (23) holds.

Since $n_1$ is order preserving and cofinal,

$$\lim_\delta \psi^{(1)}_{\phi(\varphi(\delta))}((s,t)) = \psi((s,t))\big(1 - \big(J^{1,2} + J^{2,1}\big)\big) \tag{31}$$

by (29), (27), and Lemma A.1.

By considering the two cases $J^{1,2} = 0$ and $J^{1,2} = 1$ separately, we now show that

$$\lim_\delta \psi^{(2)}_{\phi(\varphi(\delta))}((s,t)) = \psi((s,t))J^{1,2}. \tag{32}$$

If $J^{1,2} = 0$, CFs being bounded, (32) holds by (29) and (30). If $J^{1,2} = 1$, then

$$\lim_\delta \frac{n_{\phi(\varphi(\delta))}}{n_1(\phi(\varphi(\delta)))} = \lim_\delta \frac{n_2(\phi(\varphi(\delta)))}{n_1(\phi(\varphi(\delta)))} = \kappa^{-1}. \tag{33}$$

Since $n_1$ is order preserving and cofinal, using (5), (6), and Lemma A.1,

$$\lim_\delta \left( \psi_{(1)}\left( \frac{s}{\sqrt{n_1(\phi(\varphi(\delta)))}} \right) \right)^{n_1(\phi(\varphi(\delta)))} = \zeta(s),$$

implying by (33) and the continuity of exponential and logarithm,

$$\lim_\delta \left( \psi_{(1)}\left( \frac{s}{\sqrt{n_1(\phi(\varphi(\delta)))}} \right) \right)^{n_1(\phi(\varphi(\delta))) - n_{\phi(\varphi(\delta))}} = \exp\left( -\frac{(\kappa-1)s^2}{2\kappa} \right). \tag{34}$$

By (11) and Lemma A.5, $\big\{ \psi_{\xi(k)} : k \in \mathfrak{N} \big\}$ is uniformly equicontinuous; since $n_2$ is order preserving and cofinal, by (27) and Lemma A.1,

$$\lim_\delta \psi_{\xi(n_2(\phi(\varphi(\delta))))}\left( \left( e^{-\frac{1}{2}}_{\phi(\varphi(\delta))} s, t \right) \right) = \psi\left( \left( \kappa^{-\frac{1}{2}} s, t \right) \right) = \exp\left( -\frac{s^2}{2\kappa} - \frac{t^2}{2} \right). \tag{35}$$

If $J^{1,2} = 1$, (32) follows by (29), (34), and (35).

By the same argument that established (32), with the modification that $J^{2,1} = 1$ implies

$$\lim_\delta \frac{n_{\phi(\varphi(\delta))}}{n_2(\phi(\varphi(\delta)))} = \lim_\delta \frac{n_1(\phi(\varphi(\delta)))}{n_2(\phi(\varphi(\delta)))} = \kappa,$$

we obtain

$$\lim_\delta \psi^{(3)}_{\phi(\varphi(\delta))}((s,t)) = \psi((s,t))J^{2,1}. \tag{36}$$

Now (8) follows from (28), (31), (32), and (36). $\qquad \square$

<u>Proof of (9) $\Rightarrow$ (19):</u> By Lemma A.2, it suffices to show that given an arbitrary subnet $\big\{ \widehat{F}_{\phi(\beta)} : \beta \in \mathcal{F} \big\}$ of $\big\{ \widehat{F}_\alpha : \alpha \in \mathfrak{N} \times \mathfrak{N} \big\}$, there exists a further subnet $\big\{ \widehat{F}_{\phi(\varphi(\delta))} : \delta \in \mathfrak{D} \big\}$ such that

$$\lim_{\delta} \widehat{F}_{\phi(\varphi(\delta))} = \Phi. \tag{37}$$

We will establish (37) for the subnet indexed by the directed set $\mathfrak{D}$ obtained in (23).

From the definition of $v_\alpha^{(i)}$ in (21), we obtain

$$\left(v_\alpha^{(1)}\right)^2 + \left(v_\alpha^{(2)}\right)^2 = 1 \tag{38}$$

and

$$a_i = \lim_{\delta} v_{\phi(\varphi(\delta))}^{(i)} = \frac{\sigma_i}{\sqrt{\kappa^{1-i}\sigma_1^2 + \kappa^{2-i}\sigma_2^2}} \in [0,1], \tag{39}$$

with $\kappa = 0$ implying $a_1 = 1$ and $a_2 = 0$, and $\kappa = \infty$ implying $a_1 = 0$ and $a_2 = 1$. By (20), the Cauchy-Schwartz inequality, (38), and (11), $\{\widehat{F}_\alpha : \alpha \in \mathfrak{N} \times \mathfrak{N}\}$ is uniformly tight. By the LCT for nets, (37) will follow once we show, for every $u \in \Re$,

$$\lim_{\delta} \zeta_{\phi(\varphi(\delta))}(u) = \zeta(u), \tag{40}$$

where $\zeta_\alpha$ is the CF of $\widehat{F}_\alpha$. By (11) and Lemma A.5,

$$\{\psi_\alpha : \alpha \in \mathfrak{N} \times \mathfrak{N}\} \text{ is uniformly equicontinuous;} \tag{41}$$

since, by (20),

$$\zeta_\alpha(u) = \psi_\alpha\left(\left(uv_\alpha^{(1)}, \ -uv_\alpha^{(2)}\right)\right), \tag{42}$$

we obtain $\lim_{\delta} \zeta_{\phi(\varphi(\delta))}(u) = \psi((ua_1, \ -ua_2))$ from (9) (equivalently, (10)), (39), and (41), whence (40) follows from the definition of $\psi$, (39) and (38). □

<u>Proof of **(19)** ⟹ **(17)**:</u> Lemma 2 is a key step in this proof.

**Lemma 2** Assume (19) and

$$\lim_{\alpha} \widehat{G}_\alpha = \Phi, \tag{43}$$

where $\widehat{G}_\alpha \in \mathcal{M}(\Re)$ denotes the measure on $\mathcal{B}(\Re)$ induced by

$$\widehat{U}_\alpha = \frac{\left(\overline{X}_{1,n_1(\alpha)} - \mu_1\right) + \left(\overline{X}_{2,n_2(\alpha)} - \mu_2\right)}{\sqrt{\frac{\sigma_1^2}{n_1(\alpha)} + \frac{\sigma_2^2}{n_2(\alpha)}}}.$$

Then (17) holds if, with $\rho_{jk}$ as in (15),

$$\rho_{jk} = 0 \text{ for } j \neq k. \tag{44}$$

The proof of Lemma 2 is presented at the end of the section.

Since (14) $\Rightarrow$ (44), by Lemma 2 it suffices to show (19) $\Rightarrow$ (43). Let $v_\alpha^* = \left( v_\alpha^{(1)}, v_\alpha^{(2)} \right)$; then, as in (20), $\widehat{U}_\alpha = \langle Y_\alpha, v_\alpha^* \rangle$, and the argument used to establish (9) $\Rightarrow$ (19) can be repeated verbatim to conclude (9) $\Rightarrow$ (43). Since we have already proved (18) $\Rightarrow$ (9), all that remains to show is (19) $\Rightarrow$ (18).

Let $\left\{ \overline{\rho}_{\xi(k(r))} : r \in \mathfrak{N} \right\}$ be an arbitrary subsequence of $\left\{ \overline{\rho}_{\xi(k)} : k \in \mathfrak{N} \right\}$. It suffices to find a further subsequence $\left\{ \overline{\rho}_{\xi(k(r(m)))} : m \in \mathfrak{N} \right\}$ that converges to 0. By Lemma 1, there exists a further subsequence $\left\{ H_{\xi(k(r(m)))} : m \in \mathfrak{N} \right\}$ such that (24) and (25) hold. We are going to show that $\rho$ of (24) equals 0. Since $\lim_{m \in \mathfrak{N}} e_{\xi(k(r(m)))} = 1$, we obtain by (21) that

$$\lim_{m \in \mathfrak{N}} v_{\xi(k(r(m)))}^{(i)} = \sigma_i \Big/ \sqrt{\sigma_1^2 + \sigma_2^2}. \tag{45}$$

Clearly, (25) implies

$$\lim_{m \in \mathfrak{N}} \psi_{\xi(k(r(m)))}(s, t) = \psi_\rho(s, t), \tag{46}$$

where $\psi_\rho$ is the CF of $(\Phi \times \Phi)_\rho$, that is,

$$\psi_\rho(s, t) = \exp\big( -\big( s^2 + 2st\rho + t^2 \big)/2 \big). \tag{47}$$

By (42), (41), (45), and (46),

$$\lim_{m \in \mathfrak{N}} \zeta_{\xi(k(r(m)))}(u) = \exp\big( -u^2(1 - \rho^*)/2 \big), \tag{48}$$

where, with $\rho$ as in (24),

$$\rho^* = \frac{2\sigma_1 \sigma_2}{\sigma_1^2 + \sigma_2^2} \rho. \tag{49}$$

Since RHS(48) is the CF of $\mathcal{N}_{(1-\rho^*)}$ at $u$, by the LCT again,

$$\lim_{m \in \mathfrak{N}} F_{\xi(k(r(m)))} = \mathcal{N}_{(1-\rho^*)}; \tag{50}$$

since $m \mapsto \xi(k(r(m)))$ is order preserving and cofinal, by (19) and Lemma A.1, $\mathcal{N}_{(1-\rho^*)} = \Phi = \mathcal{N}_1$, showing, by (49), that $\rho = 0$, thereby establishing (18). $\qquad\square$

**Remark 4** We observed in Remark 3 that a sequence of iid copies of a random vector with dependent but uncorrelated coordinates satisfies the collection of assumptions (1), (2), (14), and (18). Clearly, if we have a sequence of iid copies of a random vector with coordinates that have a non-zero correlation coefficient, then (18) is violated and neither (9) nor (19) can hold. In this background, Proposition 1 is a generalization of the iid CLT for a sequence of random vectors.                          //

**Proposition 1** Assume (12). For all $j \geq 1$, let $0 \neq \rho = \rho_{jj}$. Then, for every subnet $\left\{ e_{\phi(\beta)} : \beta \in \mathcal{F} \right\}$ such that $\lim_{\beta} e_{\phi(\beta)} = 1$,

$$\lim_{\beta} H_{\phi(\beta)} = (\Phi \times \Phi)_{\rho} \tag{51}$$

and, with $\rho^*$ as in (49),

$$\lim_{\beta} F_{\phi(\beta)} = \mathcal{N}_{(1-\rho^*)}. \tag{52}$$

<u>Proof of Proposition 1</u> Since (12) implies (14), and (24) is satisfied by every subsequence for the same $\rho$, we conclude from Lemma 1 that

$$\lim_{k} H_{\xi(k)} = (\Phi \times \Phi)_{\rho}. \tag{53}$$

By (53) and (47),

$$\lim_{k} \psi_{\xi(k)}(s,t) = \psi_{\rho}(s,t). \tag{54}$$

Recall that the decomposition of $\psi_{\alpha}$ into the three terms in (28) was obtained solely on the basis of (14) and (5), and (5) on the basis of (1) and (2). Since (12) implies (1), (2), and (14), the decomposition of obtained in (28) continues to hold even though (18) is not true any more. With (54) substituting for (27), (31) with $\psi$ on the RHS replaced by $\psi_{\rho}$ holds. Since the limit $\kappa$ of every subnet $\{e_{\phi(\varphi(\delta))} : \delta \in \mathfrak{D}\}$ of $\{e_{\phi(\beta)} : \beta \in \mathcal{F}\}$ equals 1 by Lemma A.1, RHS(34) reduces to 1. Since (11) was obtained only by using (1) and (2), $\{\psi_{\xi(k)} : k \in \mathfrak{N}\}$ is still uniformly equicontinuous; as such, RHS(35), via (54), equals $\psi_{\rho}((s,t))$, implying that (32) with $\psi$ on the RHS replaced by $\psi_{\rho}$ holds as well. The same assertion is true for (36). Since (31), (32), and (36), with $\psi$ replaced by $\psi_{\rho}$, continue to hold, we obtain

$$\lim_{\beta} \psi_{\phi(\beta)}((s,t)) = \psi_{\rho}((s,t)),$$

implying (51) by (11) and the LCT for nets. Since $\kappa = 1$, $\lim_{\beta} v_{\phi(\beta)}^{(i)} = \sigma_i / \sqrt{\sigma_1^2 + \sigma_2^2}$ by (39); the proof of (52) follows as in the proof of (19) $\Rightarrow$ (18) above. $\qquad\square$

<u>Proof of Lemma 1</u> Compactness of $[-1,1]$ implies (24). By the LCT, (25) follows once

$$\lim_{m} H_{\xi(k(r(m)))} \circ \langle\,\cdot\,,(s,t)\rangle^{-1} = (\Phi \times \Phi)_{\rho} \circ \langle\,\cdot\,,(s,t)\rangle^{-1} \text{ for every } (s,t) \in \Re^2 \tag{55}$$

is proved by considering the cases $\tau(s,t) = 0$ and $\tau(s,t) > 0$ separately, where

$$\tau(s,t) = s^2 + 2st\rho + t^2 \text{ so that } (\Phi \times \Phi)_{\rho} \circ \langle\,\cdot\,,(s,t)\rangle^{-1} = \mathcal{N}_{\tau(s,t)}, \tag{56}$$

with $\rho$ as in (24). By (14), $\mathrm{Var}\Big(\big\langle(s,t), Y_{\xi(k)}\big\rangle\Big) = s^2 + 2st\overline{\rho}_{\xi(k)} + t^2$, implying $\tau(s,t) = \lim_{m} \mathrm{Var}\Big(\big\langle(s,t), Y_{\xi(k(r(m)))}\big\rangle\Big)$. If $\tau(s,t) = 0$, $\big\langle(s,t), Y_{\xi(k(r(m)))}\big\rangle$ converges in quadratic mean, hence in distribution, to 0, that is,

$$\lim_{m} H_{\xi(k(r(m)))} \circ \langle\,\cdot\,,(s,t)\rangle^{-1} = \text{point mass at } 0 = \mathcal{N}_0 = (\Phi \times \Phi)_{\rho} \circ \langle\,\cdot\,,(s,t)\rangle^{-1},$$

where the last equality follows from (56). If $\tau(s,t) > 0$, we will use the subsequential Lindeberg CLT of Lemma A.6 to prove (55). In the notation of Lemma A.6, with

$$W_j = \frac{s(X_{1,j} - \mu_1)}{\sigma_1} + \frac{t(X_{2,j} - \mu_2)}{\sigma_2},$$

we have

$$\sqrt{k}\,\overline{W}_k = \frac{1}{\sqrt{k}} \sum_{j=1}^{k} \left[ \frac{s(X_{1,j} - \mu_1)}{\sigma_1} + \frac{t(X_{2,j} - \mu_2)}{\sigma_2} \right] = \left\langle (s,t), Y_{\xi(k)} \right\rangle,$$

$\sigma_j^2 = s^2 + 2st\rho_{jj} + t^2$, and $\tau_k = s^2 + 2st\overline{\rho}_{\xi(k)} + t^2$. Since $\lim_m \tau_{k(r(m))} = \tau(s,t) > 0$, to apply Lemma A.6 it remains to verify that

$$\lim_m L_{k(r(m))}(\epsilon) = 0 \text{ for every } \epsilon > 0, \tag{57}$$

where, with the dependence on $(s,t)$ suppressed,

$$L_k(\epsilon) = \sum_{j=1}^{k} \mathrm{E}\left( \frac{\left[ \frac{s(X_{1,j} - \mu_1)}{\sigma_1} + \frac{t(X_{2,j} - \mu_2)}{\sigma_2} \right]^2}{k\tau_k} \left[ \frac{\left| \frac{s(X_{1,j} - \mu_1)}{\sigma_1} + \frac{t(X_{2,j} - \mu_2)}{\sigma_2} \right|}{\sqrt{k\tau_k}} > \epsilon \right] \right).$$

Recall from Lemma A.6 that $L_k(\epsilon)$ is well defined beyond a finite stage. Define

$$\mathfrak{L}_k(\epsilon, (s,t)) = \sum_{j=1}^{k} \mathrm{E}\left( \frac{\left[ \frac{s(X_{1,j} - \mu_1)}{\sigma_1} + \frac{t(X_{2,j} - \mu_2)}{\sigma_2} \right]^2}{k} \left[ \frac{\left| \frac{s(X_{1,j} - \mu_1)}{\sigma_1} + \frac{t(X_{2,j} - \mu_2)}{\sigma_2} \right|}{\sqrt{k}} > \epsilon \right] \right).$$

Note that

$$\mathfrak{L}_k \text{ is decreasing in the first argument}$$
$$\text{for } \gamma \neq 0,\ \mathfrak{L}_k(\epsilon, \gamma(s,t)) = \gamma^2 \mathfrak{L}_k\left( \frac{\epsilon}{|\gamma|}, (s,t) \right). \tag{58}$$

There exists $M^* \in \mathfrak{N}$ such that $m > M^*$ implies $\tau_{k(r(m))} > \tau(s,t)/2$. Consequently, by the definitions of $L_k$ and $\mathfrak{L}_k$, and using (58), for all $m > M^*$,

$$L_{k(r(m))}(\epsilon) < \frac{2}{\tau(s,t)} \mathfrak{L}_{k(r(m))}\left( \epsilon \sqrt{\frac{\tau(s,t)}{2}}, (s,t) \right).$$

To establish (57) from here, it suffices to show that

$$\lim_k \mathfrak{L}_k(\epsilon, (s,t)) = 0 \text{ for every } (s,t) \in \Re^2 \text{ and } \epsilon > 0. \tag{59}$$

By the identical distribution and square integrability components of (1) and (2),

$$\lim_k \mathfrak{L}_k(\epsilon, (s, t)) = 0 \text{ for every } \epsilon > 0 \tag{60}$$

holds if $(s, t)$ equals either $(1, 0)$ or $(0, 1)$. Therefore, (59) follows from (60) once we show that $L = \{(s, t) \in \Re^2 : (60) \text{ holds}\}$ is a subspace. Clearly, $L$ contains $(0, 0)$ and is closed under scalar multiplication by (58). To verify that $L$ is closed under vector addition, it suffices to show that for $(s, t), (u, v) \in \Re^2$,

$$\mathfrak{L}_k(\epsilon, (s, t) + (u, v)) \le 4\mathfrak{L}_k\left(\frac{\epsilon}{2}, (s, t)\right) + 4\mathfrak{L}_k\left(\frac{\epsilon}{2}, (u, v)\right). \tag{61}$$

Adapting an argument from the proof of Theorem 1.3 of Kundu, Majumdar, and Mukherjee (2000), since

$$\left| (s + u)\frac{(X_{1,j} - \mu_1)}{\sigma_1} + (t + v)\frac{(X_{2,j} - \mu_2)}{\sigma_2} \right|$$
$$\le 2\max\left( \left| s\frac{(X_{1,j} - \mu_1)}{\sigma_1} + t\frac{(X_{2,j} - \mu_2)}{\sigma_2} \right|, \left| u\frac{(X_{1,j} - \mu_1)}{\sigma_1} + v\frac{(X_{2,j} - \mu_2)}{\sigma_2} \right| \right),$$

$$\left( (s + u)\frac{(X_{1,j} - \mu_1)}{\sigma_1} + (t + v)\frac{(X_{2,j} - \mu_2)}{\sigma_2} \right)^2$$
$$\le 4\max\left( \left( s\frac{(X_{1,j} - \mu_1)}{\sigma_1} + t\frac{(X_{2,j} - \mu_2)}{\sigma_2} \right)^2, \left( u\frac{(X_{1,j} - \mu_1)}{\sigma_1} + v\frac{(X_{2,j} - \mu_2)}{\sigma_2} \right)^2 \right)$$
$$= 4\left[ \max\left( \left| s\frac{(X_{1,j} - \mu_1)}{\sigma_1} + t\frac{(X_{2,j} - \mu_2)}{\sigma_2} \right|, \left| u\frac{(X_{1,j} - \mu_1)}{\sigma_1} + v\frac{(X_{2,j} - \mu_2)}{\sigma_2} \right| \right) \right]^2,$$

and $(\max(x, y))^2[\max(x, y) > \epsilon] \le x^2[x > \epsilon] + y^2[y > \epsilon]$, (61) follows. By Lemma A.6, as $m \to \infty$, $\sqrt{k(r(m))}\overline{W}_{k(r(m))} = \left\langle (s, t), Y_{\xi(k(r(m)))} \right\rangle$ converges in distribution to $\mathcal{N}_{\tau(s,t)}$, and (55) follows from (56). $\square$

<u>Proof of Lemma 2</u> By Lemma A.2, given an arbitrary subnet $\{\overline{\rho}_{\phi(\beta)} : \beta \in \mathcal{F}\}$ of $\{\overline{\rho}_\alpha : \alpha \in \mathfrak{N} \times \mathfrak{N}\}$, it suffices to find a further subnet $\{\overline{\rho}_{\phi(\varphi(\delta))} : \delta \in \mathfrak{D}\}$ such that

$$\lim_\delta \overline{\rho}_{\phi(\varphi(\delta))} = 0. \tag{62}$$

Recall from (23) the existence of the directed set $\mathfrak{D}$ such that the subnet $\{e_{\phi(\varphi(\delta))} : \delta \in \mathfrak{D}\}$ converges to $\kappa$. We will show that (62) holds for $\mathfrak{D}$ by separately considering the cases $\kappa \in (0, \infty)$ and $\kappa \in \{0, \infty\}$.

By the equivalence of (i) and (iv) in the Portmanteau Theorem [van der Vaart and Wellner (1996, Theorem 1.3.4)],

$$\liminf_\alpha E\left(\widehat{W}_\alpha^2\right) \ge 1 \text{ and } \liminf_\alpha E\left(\widehat{U}_\alpha^2\right) \ge 1. \tag{63}$$

Since, by (44),

$$E\left(\widehat{W}_\alpha^2\right) = 1 - \frac{2\sigma_1\sigma_2\sqrt{e_\alpha}}{\sigma_1^2 + e_\alpha\sigma_2^2} \times \overline{\rho}_\alpha \text{ and } E\left(\widehat{U}_\alpha^2\right) = 1 + \frac{2\sigma_1\sigma_2\sqrt{e_\alpha}}{\sigma_1^2 + e_\alpha\sigma_2^2} \times \overline{\rho}_\alpha,$$

we obtain from (63)

$$\lim_\alpha \frac{2\sigma_1\sigma_2\sqrt{e_\alpha}}{\sigma_1^2 + e_\alpha\sigma_2^2} \times \overline{\rho}_\alpha = 0,$$

implying (62) for $\kappa \in (0, \infty)$. Since $|\rho_{jj}| \leq 1$ for all $j \geq 1$ implies $\left|\overline{\rho}_\alpha\right| \leq n_\alpha/m_\alpha$,

$$\left|\text{LHS}(62)\right| \leq \lim_\delta \left(n_{\phi(\varphi(\delta))}/m_{\phi(\varphi(\delta))}\right); \tag{64}$$

since $\kappa = 0$ implies RHS(64) equals the limit of $\sqrt{e_{\phi(\varphi(\delta))}}$, whereas $\kappa = \infty$ implies RHS(64) equals the limit of $1/\sqrt{e_{\phi(\varphi(\delta))}}$, (62) follows for $\kappa \in \{0, \infty\}$.                $\square$

**Remark 5** Theorem 1 stands on the Lindeberg CLT for independent random vectors. It is worth investigating if we can obtain Theorem 1 when (14) is relaxed to $m$-dependence or martingale difference array by using the CLTs under these dependence structures.        //

**Remark 6** We conjecture that the joint distribution of the standardized sample means from $k$ random samples will be asymptotically Normal in $\Re^k$ if an extension of (14) and an appropriate extension of (18) for pairwise correlation coefficients hold.        //

**A.1. Nets and subnets.** A set $\mathfrak{D}$ endowed with a reflexive, anti symmetric, and transitive binary relation $\succeq$ is called a partially ordered set. The pair $(\mathfrak{D}, \succeq)$ is called a *directed set* if, for each $\beta, \gamma \in \mathfrak{D}$, there exists $\eta \in \mathfrak{D}$ such that $\eta \succeq \beta$ and $\eta \succeq \gamma$.

Given a metric space $(S, d)$ and a directed set $(\mathfrak{D}, \succeq)$, a $S$–valued net is defined to be a function $x : \mathfrak{D} \mapsto S$; we write the net as $\{x_\beta : \beta \in \mathfrak{D}\}$. Recall that the net $\{x_\beta : \beta \in \mathfrak{D}\}$ converges to $x \in S$ if, for every $\epsilon > 0$, there exists $\beta_0(\epsilon) \in \mathfrak{D}$ such that $\beta \succeq \beta_0(\epsilon)$ implies $d(x_\beta, x) < \epsilon$. Note that a $S$–valued sequence is a $S$–valued net indexed by $\mathfrak{N}$.

Let $(\mathfrak{D}, \succeq)$ and $(\mathfrak{E}, \gg)$ be directed sets. Let $\phi : \mathfrak{E} \mapsto \mathfrak{D}$ be *order preserving*, that is, $i \gg j \Rightarrow \phi(i) \succeq \phi(j)$, and *cofinal*, that is, for each $\beta \in \mathfrak{D}$, there exists $\gamma \in \mathfrak{E}$ such that $\phi(\gamma) \succeq \beta$. Then the composite function $y = x \circ \phi$, where $x : \mathfrak{D} \mapsto S$, defines a net $\{y_\gamma : \gamma \in \mathfrak{E}\}$ in $S$, is called a subnet of $\{x_\beta : \beta \in \mathfrak{D}\}$, and is written as $\{x_{\phi(\gamma)} : \gamma \in \mathfrak{E}\}$.

**Lemma A.1** Let $\mathfrak{D}$ be a directed set and $\{x_\beta : \beta \in \mathfrak{D}\}$ a net taking values in $S$ that converges to $x \in S$. Then every subnet of $\{x_\beta : \beta \in \mathfrak{D}\}$ converges to $x$.

<u>Proof of Lemma A.1</u> This is Exercise 8 of page 188 of Munkres (2000). The proof follows from the definition of a subnet.                $\square$

**Lemma A.2** Let $(\mathfrak{D}, \succeq)$ be a directed set and $\{x_\beta : \beta \in \mathfrak{D}\}$ a net taking values in $S$. Then $\{x_\beta : \beta \in \mathfrak{D}\}$ converges to $x \in S$ if and only if every subnet of $\{x_\beta : \beta \in \mathfrak{D}\}$ has a further subnet that converges to $x$.

<u>Proof of Lemma A.2</u> Our inability to find a published proof of this very well-known result prompts us to sketch one here. The *only if* assertion follows from Lemma A.1. Conversely, suppose that every subnet of $\{x_\beta : \beta \in \mathfrak{D}\}$ has a further subnet that converges to $x$. To prove by contradiction that $\{x_\beta : \beta \in \mathfrak{D}\}$ converges to $x$, assume that $\{x_\beta : \beta \in \mathfrak{D}\}$ does not converge to $x$, that is, there exists an $\epsilon > 0$ such that for every $\beta \in \mathfrak{D}$, there exists $\beta' \in \mathfrak{D}$ satisfying $\beta' \succeq \beta$ with $d(x_{\beta'}, x) \geq \epsilon$. Note that $\mathfrak{C} = \{\beta' \in \mathfrak{D} : d(x_{\beta'}, x) \geq \epsilon\}$ is a partially ordered set with the inherited relation $\succeq$. Since $(\mathfrak{D}, \succeq)$ is a directed set, given $\beta', \gamma' \in \mathfrak{C} \subseteq \mathfrak{D}$ there exists $\eta \in \mathfrak{D}$ such that $\eta \succeq \beta'$ and $\eta \succeq \gamma'$. But given $\eta \in \mathfrak{D}$, there exists $\eta' \in \mathfrak{D}$ such that $\eta' \succeq \eta$ with $d(x_{\eta'}, x) \geq \epsilon$; that is, $\eta' \in \mathfrak{C}$. Since $\succeq$ is transitive, $\eta' \succeq \beta'$ and $\eta' \succeq \gamma'$, establishing that $(\mathfrak{C}, \succeq)$ is a directed set. Clearly, the inclusion map from $\mathfrak{C}$ to $\mathfrak{D}$ is order preserving. The argument used to establish that $(\mathfrak{C}, \succeq)$ is a directed set also establishes the cofinality of the inclusion map. Consequently $\{x_{\beta'} : \beta' \in \mathfrak{C}\}$ is a subnet of $\{x_\beta : \beta \in \mathfrak{D}\}$. By the *if* assumption, $\{x_{\beta'} : \beta' \in \mathfrak{C}\}$ has a subnet that converges to $x$, which is impossible since every element of $\{x_{\beta'} : \beta' \in \mathfrak{C}\}$ is outside an $\epsilon$-neighborhood of $x$. $\qquad\square$

**Lemma A.3** $S$ is compact if and only if every net in $S$ has a convergent subnet.

<u>Proof of Lemma A.3</u> This is the theorem stated in Exercise 10 of page 188 of Munkres (2000), who sketches an outline of the proof as a hint. $\qquad\square$

**A.2. Miscellaneous results from probability.** We have used Lemmas A.5 and A.6 in the paper; Lemma A.4 is used in the proofs of Lemmas A.5 and A.6.

**Lemma A.4** For every $\theta \in \Re$ and $m \geq 0$,

$$\left| \exp(\imath\theta) - \sum_{q=0}^{m} \frac{(\imath\theta)^q}{q!} \right| \leq \frac{|\theta|^{m+1}}{(m+1)!}. \tag{65}$$

<u>Proof of Lemma A.4</u> Again, the result is extremely well-known [Fabian and Hannan (1985, page 154)], but we could not locate a reference with a proof, hence the sketch here. The inequality in (65) is vacuously true for $\theta = 0$. Since $\exp(\imath\theta)$ and $\exp(-\imath\theta)$, and, for every $q \geq 1$, $(\imath)^q$ and $(-\imath)^q$, are complex conjugates of each other, it suffices to prove (65) for $\theta > 0$. With $c$ denoting the cosine function and $s$ the sine function,

$$\exp(\imath\theta) - \sum_{q=0}^{m} \frac{(\imath\theta)^q}{q!} = c(\theta) - \sum_{q=0}^{m} \frac{\theta^q}{q!} c^{(q)}(0) + \imath \left[ s(\theta) - \sum_{q=0}^{m} \frac{\theta^q}{q!} s^{(q)}(0) \right];$$

that is because for $m \geq 0$,

$$c^{(q)} = \begin{cases} c & \text{if } q = 4m \\ -s & \text{if } q = 4m+1 \\ -c & \text{if } q = 4m+2 \\ s & \text{if } q = 4m+3 \end{cases} \quad \text{and} \quad s^{(q)} = \begin{cases} s & \text{if } q = 4m \\ c & \text{if } q = 4m+1 \\ -s & \text{if } q = 4m+2 \\ -c & \text{if } q = 4m+3. \end{cases}$$

By Theorem 7.6 of Apostol (1967), with $g$ equal to either $c$ or $s$,

$$g(\theta) - \sum_{q=0}^{m} \frac{\theta^q}{q!} g^{(q)}(0) = \frac{1}{m!} \int_0^\theta (\theta - t)^m g^{(m+1)}(t) dt.$$

Thus, substituting $u = \theta - t$,

$$\text{LHS}(65) = \frac{1}{m!} \left[ \left( \int_0^\theta u^m c^{(m+1)}(\theta - u) du \right)^2 + \left( \int_0^\theta u^m s^{(m+1)}(\theta - u) du \right)^2 \right]^{1/2}.$$

Consider the probability distribution $\mu$ on $[0, \theta]$ with Lebesgue density $u^m(m+1)/\theta^{m+1}$; then (65) follows from the variance inequality with respect to $\mu$ and the fact that for every $k \geq 0$ and every $x \in \Re$, $\left(c^{(k)}(x)\right)^2 + \left(s^{(k)}(x)\right)^2 = 1$. □

**Lemma A.5** Let $\{\nu_g : g \in \mathcal{G}\} \subset \mathcal{M}\left(\Re^k\right)$ be uniformly tight. Let $\psi_g$ be the CF of $\nu_g$. Then $\{\psi_g : g \in \mathcal{G}\}$ is a uniformly equicontinuous family of functions.

<u>Proof of Lemma A.5</u> Fix $\epsilon > 0$ arbitrarily. Let $K$ be a compact subset of $\Re^k$ such that

$$\nu_g(K^c) < \frac{\epsilon}{4} \text{ for all } g \in \mathcal{G}. \tag{66}$$

By Lemma A.4 (with $m = 0$) and Cauchy-Schwartz inequality, for all $s, t, x \in \Re^k$

$$\left| 1 - \exp\left(\imath \langle t - s, x \rangle\right) \right| \leq \left| \langle t - s, x \rangle \right| \leq \|t - s\| \|x\|.$$

There exists $\delta = \epsilon/2M > 0$ such that $\|t - s\| < \delta$ implies

$$\sup\left\{ \left| 1 - \exp\left(\imath \langle t - s, x \rangle\right) \right| : x \in K \right\} < \frac{\epsilon}{2}, \tag{67}$$

where $M$ bounds the compact set $K$. Therefore, by (66) and (67), $\|s - t\| < \delta$ implies

$$\left| \psi_g(s) - \psi_g(t) \right| < \frac{\epsilon}{2} \nu_g(K) + 2\nu_g(K^c) < \epsilon,$$

completing the proof of the lemma. □

**Lemma A.6** (*A subsequential Lindeberg CLT*) Let $\{W_j : j \in \mathfrak{N}\}$ be a sequence of independent random variables with $\text{E}(W_j) = 0$ and $\text{Var}(W_j) = \sigma_j^2$. For $k \in \mathfrak{N}$, let

$$\overline{W}_k = \frac{1}{k}\sum_{j=1}^{k} W_j \text{ and } \tau_k = \frac{1}{k}\sum_{j=1}^{k}\sigma_j^2 = \text{Var}\left(\sqrt{k}\overline{W}_k\right).$$

Let $\left\{\tau_{k(r)} : r \in \mathfrak{N}\right\}$ be a subsequence of $\{\tau_k : k \in \mathfrak{N}\}$ such that

$$\lim_r \tau_{k(r)} = \tau > 0. \tag{68}$$

Then,

$$\text{there exists } K^* \in \mathfrak{N} \text{ such that } \tau_k > 0 \text{ for all } k \geq K^*. \tag{69}$$

For $\epsilon > 0$ and $k \geq K^*$, let

$$L_k(\epsilon) = \sum_{j=1}^{k} \text{E}\left(\left((k\tau_k)^{-\frac{1}{2}}W_j\right)^2\left[\left|(k\tau_k)^{-\frac{1}{2}}W_j\right| > \epsilon\right]\right).$$

Assume

$$\lim_r L_{k(r)}(\epsilon) = 0 \text{ for every } \epsilon > 0. \tag{70}$$

Then, as $r \to \infty$,

$$\sqrt{k(r)}\overline{W}_{k(r)} \text{ converges in distribution to } \mathcal{N}_\tau. \tag{71}$$

<u>Proof of Lemma A.6</u> First note that (68) rules out the possibility that $\tau_k = 0 \ \forall \ k \in \mathfrak{N}$. Thus, there exists $K^* \in \mathfrak{N}$ such that $\tau_{K^*} > 0$. Since $\{k\tau_k : k \in \mathfrak{N}\}$ is a nondecreasing sequence, (69) follows. Thus we can, without loss of generality, assume that $K^* = 1$.

For every $k \in \mathfrak{N}$ and $1 \leq j \leq k$, let $\{Z_{k,j} : 1 \leq j \leq k\}$ be iid $\sim \Phi$; further, let

$$Y_{k,j} = (k\tau_k)^{-\frac{1}{2}}W_j, \ \phi_{k,j} = \text{the CF of } Y_{k,j}, \ S_k = \sum_{j=1}^{k} Y_{k,j}, \text{ and } \Upsilon_k = \text{the CF of } S_k$$

$$\varphi_{k,j}^2 = \text{Var}(Y_{k,j}) = \text{E}\left(Y_{k,j}^2\right) = (k\tau_k)^{-1}\sigma_j^2 \text{ so that } \sum_{j=1}^{k}\varphi_{k,j}^2 = 1 \tag{72}$$

$$\widehat{Y}_{k,j} = \varphi_{k,j}Z_{k,j}, \ \widehat{\phi}_{k,j} = \text{the CF of } \widehat{Y}_{k,j}, \ \widehat{S}_k = \sum_{j=1}^{k}\widehat{Y}_{k,j}, \text{ and } \widehat{\Upsilon}_k = \text{the CF of } \widehat{S}_k.$$

Clearly, for every $t \in \Re$,

$$\Upsilon_k(t) = \prod_{j=1}^{k}\phi_{k,j}(t), \ \widehat{\phi}_{k,j}(t) = \exp\left(-\frac{1}{2}t^2\varphi_{k,j}^2\right), \text{ and } \widehat{\Upsilon}_k(t) = \exp\left(-\frac{1}{2}t^2\right).$$

Fix $0 \neq t \in \Re$ arbitrarily. We will show that

$$\lim_{T} \left| \Upsilon_{k(r)}(t) - \prod_{j=1}^{k(r)} \psi_{k(r),j}(t) \right| = 0 = \lim_{T} \left| \widehat{\Upsilon}_{k(r)}(t) - \prod_{j=1}^{k(r)} \psi_{k(r),j}(t) \right|, \tag{73}$$

where

$$\psi_{k,j}(t) = 1 - \frac{1}{2} t^2 \varphi_{k,j}^2, \tag{74}$$

which implies

$$\lim_{T} \Upsilon_{k(r)}(t) = \exp\left( -\frac{1}{2} t^2 \right). \tag{75}$$

Note that (75) implies, by the LCT, Proposition 9.3.4 of Dudley (1989), and Lemma A.5, that $\{ \Upsilon_{k(r)} : r \in \mathfrak{N} \}$ is uniformly equicontinuous, whence, by (68),

$$\lim_{T} \Upsilon_{k(r)}\left( t \sqrt{\tau_{k(r)}} \right) = \exp\left( -\frac{1}{2} t^2 \tau \right),$$

implying, since $\sqrt{k} \overline{W}_k = \sqrt{\tau_k} S_k$, (71) by the LCT. Thus, it remains to show (73).

Since $\{ Y_{k,j} : 1 \le j \le k \}$ is a row in the sense of Definition 4.1.4 of Fabian and Hannan (1985), for arbitrary $\epsilon > 0$, by Lemma 4.1.7 of Fabian and Hannan (1985),

$$a_k^2 = \max\{ \varphi_{k,j}^2 : 1 \le j \le k \} \le \epsilon^2 + L_k(\epsilon), \tag{76}$$

implying, by (70),

$$\lim_{T} a_{k(r)} = 0. \tag{77}$$

By (72) and (76), $\widehat{L}_k(\epsilon) = \sum_{j=1}^{k} \mathrm{E}\left( \widehat{Y}_{k,j}^2 \left[ \left| \widehat{Y}_{k,j} \right| > \epsilon \right] \right) \le E\left( Z^2 \left[ \left| Z \right| > \frac{\epsilon}{a_k} \right] \right)$; by (77),

$$\lim_{T} \widehat{L}_{k(r)}(\epsilon) = 0 \text{ for every } \epsilon > 0. \tag{78}$$

With $h(y) = \exp(\imath t y) - 1 - \imath t y + t^2 y^2 / 2$ for $y \in \Re$, we obtain

$$\mathrm{E}(h(Y_{k,j})) = \phi_{k,j}(t) - \psi_{k,j}(t)$$
$$\mathrm{E}\left( h\left( \widehat{Y}_{k,j} \right) \right) = \widehat{\phi}_{k,j}(t) - \psi_{k,j}(t).$$

By (65) with $m = 1$, $\left| h(y) \right| \le t^2 y^2$; with $m = 2$, $\left| h(y) \right| \le \left| t \right|^3 |y| y^2$. For every $\epsilon > 0$,

$$\sum_{j=1}^{k} \left| \phi_{k,j}(t) - \psi_{k,j}(t) \right| \le \sum_{j=1}^{k} \left[ \epsilon \left| t \right|^3 \varphi_{k,j}^2 + t^2 \mathrm{E}\left( Y_{k,j}^2 \left[ \left| Y_{k,j} \right| > \epsilon \right] \right) \right] = \epsilon \left| t \right|^3 + t^2 L_k(\epsilon);$$

similarly,

$$\sum_{j=1}^{k} \left| \widehat{\phi}_{k,j}(t) - \psi_{k,j}(t) \right| \le \epsilon |t|^3 + t^2 \widehat{L}_k(\epsilon).$$

By (77), there exists $R^* = R^*(t) \in \mathfrak{N}$ such that $r > R^*$ implies, for all $1 \le j \le k(r)$,

$$0 \le \varphi_{k(r),j}^2 \le \frac{4}{t^2},$$

so that by (74), for $r > R^*$, $\left| \psi_{k(r),j}(t) \right| \le 1$. By Lemma 4.1.2 of Fabian and Hannan (1985), for $r > R^*$,

$$\left| \prod_{j=1}^{k(r)} \phi_{k(r),j}(t) - \prod_{j=1}^{k(r)} \psi_{k(r),j}(t) \right| \le \sum_{j=1}^{k} \left| \phi_{k,j}(t) - \psi_{k,j}(t) \right| \le \epsilon |t|^3 + t^2 L_{k(r)}(\epsilon)$$

$$\left| \prod_{j=1}^{k(r)} \widehat{\phi}_{k(r),j}(t) - \prod_{j=1}^{k(r)} \psi_{k(r),j}(t) \right| \le \sum_{j=1}^{k} \left| \widehat{\phi}_{k,j}(t) - \psi_{k,j}(t) \right| \le \epsilon |t|^3 + t^2 \widehat{L}_{k(r)}(\epsilon);$$

that establishes (73) by (70) and (78), completing the proof. $\qquad\square$

RAJESHWARI MAJUMDAR                        SUMAN MAJUMDAR
rajeshwari.majumdar@uconn.edu              suman.majumdar@uconn.edu
PO Box 47                                  1 University Place
Coventry, CT 06238                         Stamford, CT 06901